\documentclass[
aps,%
12pt,%
final,%
notitlepage,%
oneside,%
onecolumn,%
nobibnotes,%
nofootinbib,%
superscriptaddress,%
showpacs,%
centertags]%
{revtex4}

\def\Boxx{
\vbox{
\halign to5.8pt 
{\strut##&
\hfil ## \hfil \cr
&$\kern -0.5pt
\sqcap$ \cr
\noalign{\kern -5pt
\hrule}
}}}

\def\bea{\begin{eqnarray}}
\def\eqnn{\bea\label}
\def\eea{\end{eqnarray}}
\def\nn{\nonumber}

\newcommand{\eqna}[1]{\begin{subequations} \label{#1}
\begin{eqnarray}}
\def\eena{\end{eqnarray}
\end{subequations}}

\def\hgamma{\hat \gamma}
\def\tgamma{\tilde \gamma}

 \def\y{\eta}
 \def\s{\sigma}  \def\t{\tau}
 
\def\ra{\rightarrow}
\def\hvf{\hat{\varphi}}  
\def\tvf{\tilde{\varphi}}  
\def\bbz{Z\!\!\!Z}
\def\k{\kappa}
\def\Uf{U_q(sl(4))}
\def\hd{{\hat{\cal D}}}
\def\({\left(}
\def\){\right)}
\def\bbn{I\!\!N}
\def\bz{{\bar z}}
\def\G{\Gamma}

\def\cl{{\cal L}}
\def\be{\begin{equation}}
\def\eqn{\begin{equation}\label}
\def\ee{\end{equation}}

\def\bea{\begin{eqnarray}}
\def\eqnn{\bea\label}
\def\eea{\end{eqnarray}}
\def\nn{\nonumber}

\def\half{{\textstyle{1\over2}}}
\def\ha{{\textstyle{1\over2}}}
\def\sixth{{\textstyle{1\over6}}}
\def\bv{{\bar v}}
\def\hm{{\hat M}}

\def\vf{\varphi} 
\def\pd{\partial} 
\def\om{\omega}  \def\hh{{\hat h}}
\def\l{\lambda}

\begin{document}
\title{q-Plane Wave Solutions of q-Weyl gravity}
\author{\firstname{Vladimir K.}~\surname{Dobrev}}
\email{dobrev@inrne.bas.bg} \affiliation{Institute of Nuclear
Research and Nuclear Energy, Bulgarian Academy of Sciences, 72
Tsarigradsko Chaussee, 1784 Sofia, Bulgaria\\
and\\
Abdus Salam International Center for Theoretical Physics,\\ Strada
Costiera 11, 34100 Trieste, Italy}
\author{\firstname{Stephen G.}~\surname{Mihov}}
\email{smikhov@inrne.bas.bg} \affiliation{Institute of Nuclear
Research and Nuclear Energy, Bulgarian Academy of Sciences, 72
Tsarigradsko Chaussee, 1784 Sofia, Bulgaria}

  \begin{abstract}
We give  solutions of the q-deformed equations of quantum conformal
Weyl gravity in terms of q-deformed plane waves.
\end{abstract}

 \maketitle

\section{Introduction}

One of the purposes of quantum deformations
is to provide an alternative of the regularization
procedures of quantum field theory. Applied to
Minkowski space-time the quantum deformations approach
is also an alternative to Connes' noncommutative geometry \cite{Con}.
The first step in such an approach is to construct
a noncommutative quantum deformation of Minkowski space-time.
There are several possible such deformations,
cf. \cite{CSSW,SWZ,Mja,Mjb,VKD2}.
We shall follow the deformation of \cite{VKD2}
which is different from the others, the most
important aspect being that it is related
to a deformation of the conformal group.

The first problem to tackle in a noncommutative deformed setting is
to study the $q$-deformed analogues of the conformally invariant
equations. Here we continue the study of hierarchies of deformed
equations derived in \cite{VKD2,VKD3,VKD4} with the use of quantum
conformal symmetry. We give now a description of our setting
starting from the simplest example.

It is well known that the d'Alembert equation
\eqn{dal} \Boxx ~\vf(x) ~~=~~ 0 ~, \qquad \Boxx ~=~ \pd^\mu \pd_\mu
~=~ (\pd_0)^2 - (\vec\pd)^2 ~, \end{equation}
is conformally invariant,
cf., e.g., \cite{BR}. Here ~$\vf$~ is a scalar field of fixed conformal
weight, ~$x\ =\ (x_0,x_1,x_2,x_3)$~ denotes the Minkowski
space-time coordinates. Not known was the fact that (\ref{dal})\
may be interpreted as conditionally conformally invariant
equation and thus may be rederived from a subsingular vector of a
Verma module of the algebra\ $sl(4)$, the complexification of the
conformal algebra\ $su(2,2)$ \cite{VKD3}.

The same idea was used in \cite{VKD3} to derive a $q$-d'Alembert
equation, namely, as arising from a subsingular vector of a Verma
module of the quantum algebra\ $\Uf$. The resulting equation is a
$q$-difference equation and the solution spaces are built on the
noncommutative $q$-Minkowski space-time of \cite{VKD2}.

Besides the $q$-d'Alembert equation in \cite{VKD3} were derived a
whole hierarchy of equations corresponding to the massless
representations of the conformal group and parametrized by a
nonnegative integer \ $r$\ \cite{VKD3}. The case $r=0$ corresponds
to the $q$-d'Alembert equation, while for each \ $r>0$ \ there are
two couples of equations involving fields of conjugated Lorentz
representations of dimension $r+1$. For instance, the case ~$r=1$~
corresponds to the massless Dirac equation, one couple of equations
describing the neutrino, the other couple of equations describing
the antineutrino, while the case ~$r=2$~ corresponds to the Maxwell
equations.

The construction of solutions of the $q$-d'Alembert hierarchy was
started in
 \cite{DK} with the $q$-d'Alembert equation. One of the solutions given
was a deformation of the plane wave as a formal power series in the
noncommutative coordinates of $q$-Minkowski space-time and
four-momenta. This $q$-plane wave has some properties analogous to
the classical one but is not an exponent or $q$-exponent. Thus, it
differs conceptually from the classical plane wave and may serve as
a regularization of the latter. In the same sense it differs from
the $q$-plane wave in the paper \cite{Mey}, which is not surprising,
since there is used different $q$-Minkowski space-time (from
\cite{CSSW,SWZ,Mja} and different $q$-d'Alembert equation both based
only on a (different) $q$-Lorentz algebra, and not on $q$-conformal
(or $\Uf$) symmetry as in our case. In fact, it is not clear whether
the $q$-Lorentz algebra of \cite{CSSW,SWZ,Mja} used in \cite{Mey} is
extendable to a $q$-conformal algebra.

For the equations labelled by $r>0$ it turned out that one needs a
second $q$-deformation of the plane wave in a conjugated basis
 \cite{DGPZ}. The solutions of the hierarchy in terms of the two
$q$-plane waves were given in \cite{DGPZ} for $r=1$ and in
 \cite{DZ} for $r>1$. Later these two $q$-plane waves were generalized
and correspondingly more general solutions of the hierarchy were
given in \cite{DPZa}.

Another hierarchy derived in \cite{VKD2} is the Maxwell hierarchy.
The two hierarchies have only one common member - the Maxwell
equations - they are the lowest member of the Maxwell hierarchy and
the $r=2$ member of the massless hierarchy. The compatibility of the
solutions of the free $q$-Maxwell equations with the $q$-potential
equations was studied in \cite{DPZb}.

Another family contained in \cite{VKD4}, but not explicated there,
is related to the linear conformal  Weyl  gravity. Its study started
in \cite{DP}, where was written down the quantum conformal
deformation of the linear conformal  Weyl  gravity. In the present
paper we continue this study by constructing solutions of these
$q$-deformed equations.

\section{Linear conformal gravity}

We shall consider the quantum group  analogs of linear conformal
gravity following the approach of \cite{VKD4}. We start with the
$q=1$ situation and we first write the Weyl gravity equations in an
indexless formulation, trading the indices for two conjugate
variables $z,\bz$.

Weyl gravity is governed by the Weyl tensor: \eqn{wt} C_{\mu\nu\s\t}
= R_{\mu\nu\s\t} -\half( g_{\mu\s}R_{\nu\t}+ g_{\nu\t}R_{\mu\s}-
g_{\mu\t}R_{\nu\s}- g_{\nu\s}R_{\mu\t}) +\sixth (g_{\mu\s}g_{\nu\t}-
g_{\mu\t}g_{\nu\s})R \ ,\ee where $g_{\mu\nu}$ is the metric tensor.
Linear conformal gravity is obtained when the metric tensor is
written as: $g_{\mu\nu} = \y_{\mu\nu} + h_{\mu\nu}$, where
$\y_{\mu\nu}$ is the flat Minkowski metric, $h_{\mu\nu}$ are small
so that all quadratic and higher order terms are neglected. In
particular:  $R_{\mu\nu\s\t} = \half(\pd_\mu \pd_\t h_{\nu\s}+
\pd_\nu \pd_\s h_{\mu\t}- \pd_\mu \pd_\s h_{\nu\t}- \pd_\nu \pd_\t
h_{\mu\s})$. The equations of linear conformal gravity are:
\eqn{wfr}  \pd^\nu \pd^\t C_{\mu\nu\s\t} = T_{\mu\s}\ , \ee where
$T_{\mu\nu}$ is the energy-momentum tensor. {}From the symmetry
properties of the Weyl tensor it follows that it has ten independent
components. These may be chosen as follows (introducing notation for
future use): \eqnn{weyli} &C_0=C_{0123}\ , \quad C_1=C_{2121}\
,\quad C_2=C_{0202}\ ,\quad
C_3=C_{3012}\ ,\nn\\
&C_4=C_{2021}\ ,\quad C_5=C_{1012}\ ,\quad
C_6=C_{2023}\ ,\nn\\
&C_7=C_{3132}\ ,\quad C_8=C_{2123}\ ,\quad C_9=C_{1213}\ . \eea
Furthermore, the Weyl tensor  transforms as the direct sum of two
conjugate Lorentz irreps, which we shall denote as $C^\pm$. The
tensors $T_{\mu\nu}$ and $h_{\mu\nu}$ are symmetric and traceless
with nine independent components.

In order to be more precise we recall that the physically relevant
representations $T^\chi$ of the 4-dimensional conformal algebra
$su(2,2)$ may be labelled by $\chi = [n_1,n_2;d]$, where $n_1, n_2$
are non-negative integers fixing finite-dimensional irreducible
representations of the Lorentz subalgebra, (the dimension being
$(n_1 +1)(n_2 +1)$), and ~$d$~ is the conformal dimension (or
energy). (In the literature these Lorentz representations are
labelled also by $(j_1,j_2)=(n_1/2,n_2/2)$.) The Weyl tensor
transforms as the direct sum: \eqnn{sgw} &\chi^+ \oplus \chi^- \cr
&\cr &\chi^+ = [4,0;2] ~, \quad    \chi^- = [0,4;2] ~, \eea while
the energy-momentum tensor and the metric transform as: \eqn{sgt}
\chi_T = [2,2;4] ~, \qquad \chi_h = [2,2;0] ~, \ee as anticipated.
Indeed, $(n_1,n_2) = (2,2)$ is the nine-dimensional Lorentz
representation, (carried by $T_{\mu\nu}$ or $h_{\mu\nu}$), and
$(n_1,n_2) = (4,0),(0,4)$ are the two conjugate five-dimensional
Lorentz representations, (carried by $C^\pm$), while the conformal
dimensions are the canonical dimensions of a energy-momentum tensor
($d=4$), of the metric ($d=0$), and of the Weyl tensor ($d=2$).

As we mentioned in the Introduction the case of Weyl gravity belongs
together with the Maxwell case to an infinite family parametrized by
~$n=1,2,\dots$~ where the signatures analogous to (\ref{sgw}),
(\ref{sgt}) are: \eqnn{sgww} &\chi^+_m = [2m,0;2] ~, \quad    \chi^- =
[0,2m;2] ~, \\ &\chi^m_T = [m,m;2+m] ~, \qquad \chi^m_h =
[m,m;2-m] ~. \eea The Maxwell case is obtained for $m=1$, Weyl
gravity for $m=2$. (Note, however, that the representations
$\chi^m_h$ are not unitary for $m>2$.)

Further, we shall use  the fact that a Lorentz irrep (spin-tensor)
with signature $(n_1,n_2)$ may be represented by a polynomial
$G(z,\bz)$ in $z,\bz$ of order $n_1,n_2$, resp. More explicitly, for
the Weyl gravity representations mentioned above we use \cite{DP}:
\eqnn{deff}
C^+(z) &=& z^4C^+_4+z^3C^+_3+z^2C^+_2+zC^+_1+C^+_0 \ ,\\
C^-(\bz) &=& \bz^4C^-_4+\bz^3C^-_3+\bz^2C^-_2+\bz C^-_1+C^-_0 \ ,\nn\\
T(z,\bz )& =& z^2\bz^2T'_{22}+z^2\bz T'_{21}+z^2T'_{20}+\cr &&
+z\bz^2T'_{12}+z\bz T'_{11}+zT'_{10}+ \cr
&& +\bz^2T'_{02}+\bz T'_{01}+T'_{00} \ ,\\
h(z,\bz ) &=& z^2\bz^2h'_{22}+z^2\bz h'_{21}+z^2h'_{20}+\cr &&
+z\bz^2h'_{12}+z\bz h'_{11}+zh'_{10}+\cr && +\bz^2h'_{02}+\bz
h'_{01}+h'_{00} \ , \eea where the indices on the RHS are not
Lorentz-covariance indices, they just indicate the powers of
$z,\bz$. The components $C^\pm_k$ are given in terms of the Weyl
tensor components as follows \cite{DP}: \eqnn{wcomp} &&C^+_0=C_2-\ha
C_1-C_6+i(C_0+ \ha C_3+ C_7)\cr &&C^+_1=2(C_4-C_8+i(C_9-C_5))\cr
&&C^+_2=3(C_1-iC_3)\cr &&C^+_3=8(C_4+C_8+i(C_9+C_5))\cr &&C^+_4=C_2-
\ha C_1+C_6+i(C_0+ \ha C_3 - C_7)\cr &&C^-_0=C_2- \ha C_1
-C_6-i(C_0+ \ha C_3+ C_7)\cr &&C^-_1=2(C_4-C_8-i(C_9-C_5))\cr
&&C^-_2=3(C_1+iC_3)\cr &&C^-_3=2(C_4+C_8-i(C_9+C_5))\cr &&C^-_4=C_2-
\ha C_1 +C_6-i(C_0+ \ha C_3- C_7) \eea while the components
$T'_{ij}$ are given in terms of $T_{\mu\nu}$ as follows \cite{DP}:
\eqnn{enem} &&T'_{22} = T_{00} + 2 T_{03} + T_{33} \cr &&T'_{11} =
T_{00} - T_{33} \cr &&T'_{00} = T_{00} - 2 T_{03} + T_{33} \cr
&&T'_{21} = T_{01} + iT_{02} + T_{13} + iT_{23} \cr &&T'_{12} =
T_{01} - iT_{02} + T_{13} - iT_{23} \cr &&T'_{10} = T_{01} + iT_{02}
- T_{13} - iT_{23} \cr &&T'_{01} = T_{01} - iT_{02} - T_{13} +
iT_{23} \cr &&T'_{20} = T_{11} + 2iT_{12} - T_{22} \cr &&T'_{02} =
T_{11} - 2iT_{12} - T_{22} \eea and similarly for $h'_{ij}$ in terms
of $h_{\mu\nu}\,$.

In these terms all linear conformal  Weyl  gravity equations
(\ref{wfr}) may be written in compact form as the following pair of
equations: \eqn{wge} I^+ ~ C^+(z) ~=~ T(z,\bz) ~, \qquad I^- ~
C^-(\bz) ~=~ T(z,\bz) ~, \ee where  the operators $I^\pm$ are given
as follows: \eqnn{oper} I^+ &=&
\Bigl(z^2\bz^2\pd_+^2+z^2\pd^2_v+\bz^2\pd_{\bv }^2+\pd^2_-+\cr &&
+2z^2\bz \pd_v\pd_++2z\bz^2\pd_+\pd_{\bv }+2z\bz
(\pd_-\pd_++\pd_v\pd_{\bv })+\cr && +2\bz \pd_-\pd_{\bv
}+2z\pd_v\pd_-\Bigr)\pd^2_z -\cr &&
-6\Bigl(z\bz^2\pd_+^2+z\pd^2_v+2z\bz \pd_v\pd_++\bz^2\pd_+\pd_{\bv
}+\cr && +\bz(\pd_-\pd_++\pd_v\pd_{\bv })+\pd_v\pd_-\Bigr)\pd_z+\cr
&& 12\Bigl(\bz^2\pd_+^2+\pd_v^2+2\bz\pd_v\pd_+\Bigr)\ , \\
I^{-} &=& \Bigl(z^2\bz^2\pd_+^2+z^2\pd^2_v+\bz^2\pd_{\bv
}^2+\pd^2_-+\cr && +2z^2\bz \pd_v\pd_++2z\bz^2\pd_+\pd_{\bv }+2z\bz
(\pd_-\pd_++\pd_v\pd_{\bv })+\cr && +2\bz \pd_-\pd_{\bv
}+2z\pd_v\pd_-\Bigr)\pd^2_{\bz } -\cr && -6\Bigl(z^2\bz
\pd_+^2+\bz\pd^2_{\bv }+2z\bz \pd_+\pd_{\bv}+ z^2\pd_v\pd_++\cr &&
+z(\pd_-\pd_++\pd_v\pd_{\bv })+\pd_-\pd_{\bv }\Bigr)\pd_{\bz }+ \cr
&& 12\Bigl(z^2\pd_+^2+\pd_{\bv }^2+2z\pd_+\pd_{\bv }\Bigr) \ \nn
\eea where the variables ~$x_\pm,v,\bv$~ are expressed through the
Minkowski coordinates ~$x_0,x_1,x_2,x_3$~ as follows \cite{VKD2}:
\eqn{mink} x_\pm \ \equiv \ x_0 \pm x_3\ , \qquad v \ \equiv \ x_1
-i x_2\ , \qquad \bv \ \equiv \ x_1 + i x_2\ . \ee These variables
have, (unlike the $x_\mu$), definite group--theoretical
interpretation as part of a six-dimensional coset of the conformal
group\ $SU(2,2)$ (as explained in \cite{VKD2}). In terms of these
variables, e.g., the d'Alembert equation (\ref{dal}) is: \eqn{dali}
\Boxx\ \vf = (\pd_- \pd_+ \ -\ \pd_{v}\ \pd_\bv )\ \vf \ =\ 0 \ .\ee

To make more transparent the origin of (\ref{wge}) and in the same
time to derive the quantum group deformation of (\ref{wge}),
(\ref{oper}) we first introduce the following parameter-dependent
operators: \eqnn{opsb} I^+(n) ~&=&~ \ha \Bigl( n(n-1) I^2_1 I^2_2
- 2(n^2-1) I_1 I^2_2 I_1  + n(n+1) I^2_2 I^2_1 \Bigr)
~, \\
I^-(n) ~&=&~ \ha \Bigl( n(n-1) I^2_3 I^2_2 - 2(n^2-1) I_3 I^2_2 I_3
+ n(n+1) I^2_2 I^2_3 \Bigr)~,\nn \eea where \eqn{opsbb}
 I_1 ~\equiv ~\pd_z ~, \quad I_2 ~\equiv ~ \bz z \pd_+ + z\pd_v
+ \bz \pd_\bv + \pd_-  ~, \quad I_3 ~\equiv ~\pd_\bz ~. \ee It is
easy to check that we have the following relation: \eqn{rel} I^\pm =
I^\pm(4) \ , \ee i.e., (\ref{wge}) are written as: \eqn{wgee} I^+(4)
\ C^+(z) ~=~ T(z,\bz) ~, \qquad I^-(4)\  C^-(\bz) ~=~ T(z,\bz) ~.
\ee

 We note in passing that group-theoretically the operators $I_a$
correspond to the three simple roots of the root system of $sl(4)$,
while the operators $I^\pm_n$ correspond to the two non-simple
non-highest roots \cite{VKDcl}.

This is the form that is immediately generalizable to the
$q$-deformed case. We first present the necessary formalism in the
next Section.

Using the same operators we can write down the pair of equations which give the Weyl
tensor components in terms of the metric tensor:
\eqn{wme} I^+(2) ~ h(z,\bz) ~=~ C^-(\bz) ~, \qquad
I^-(2) ~ h(z,\bz) ~=~ C^+(z) ~. \ee

We stress the advantage of the indexless formalism due to which two
different pairs of equations, (\ref{wgee}), (\ref{wme}), may be
written using the same parameter-dependent operator expressions by
just specializing the values of the parameter.

The analogues of (\ref{wge}), (\ref{wme}), for the family
(\ref{sgww}) are: \eqnn{wgem} I^+_m(4) ~ C^+_m(z) ~=~ T_m(z,\bz) ~,
\qquad I^-_m(4) ~
C^-_m(\bz) ~=~ T_m(z,\bz) ~, \\
I^+_m(2) ~ h_m(z,\bz) ~=~ C^-_m(\bz) ~, \qquad
I^-_m(2) ~ h_m(z,\bz) ~=~ C^+_m(z) ~, \eea
where the operators ~$I^+_m(n)$~ are of order ~$m$, and can be
found in \cite{Dosv} also in the $q$-deformed case.

\section{$q$-deformed setting}

In the $q$-deformed case we use the noncommutative $q$-Minkowski
space-time of
 \cite{VKD2} which is
given by the following commutation relations (with $\l \equiv
q-q^{-1}$): \eqn{coop} x_\pm v = q^{\pm 1} v x_\pm \,, \ \
 x_\pm \bv = q^{\pm 1} \bv x_\pm \,,\ \
x_+ x_- - x_- x_+ = \l v\bv \,, \ \ \bv v = v\bv \,,
\ee
with the deformation parameter being a phase:\ $\vert q\vert =1$.
Relations (\ref{coop}) are preserved by the anti-linear anti-involution
 \ $\om$ \ :
\eqn{cnjm} \om (x_\pm) \ = \ x_\pm \ , \quad \om (v) \ = \ \bv \ ,
\quad \om (q) \ = \ \bar q \ = \ q^{-1} \ ,
\quad (\om (\l) \ = \ -\l ) \ . \ee

The solution spaces consist of formal power series in the
$q$-Minkowski coordinates (which we give in two conjugate bases):
\bea &&\vf \ = \ \sum_{j,n,\ell,m\in\bbz_+}\ \mu_{j n\ell m}\ \vf_{j
n\ell m} \ , \qquad \vf_{j n\ell m}\ =\
\hvf_{j n\ell m},\ \tvf_{j n\ell m} \ ,\label{spc}\\
&&\hvf_{j n\ell m} \ = \ v^j\ x^n_-\ x^\ell_+\ \bv^m \ ,\label{spca}\\ [2mm]
&&\tvf_{j n\ell m} \ = \ \bv^m\ x^\ell_+\ x^n_-\ v^j \ = \
\om (\hvf_{j n\ell m}) \ . \label{spcb}\eea
The solution spaces (\ref{spc}) are representation spaces of the
quantum algebra $U_q(sl(4))$. For the latter we use the rational
basis of Jimbo \cite{MJ}.
The action of\ $\Uf$\ on\ $\hvf_{j n\ell m}$\ was given in
 \cite{VKD1}, and on \ $\tvf_{j n\ell m}$ in \cite{DGPZ}.
Because of the conjugation $\om$ we are actually working with
the conformal quantum algebra which is a deformation of $U(su(2,2))$.

 Further we suppose that\ $q$\ is not a nontrivial root of unity.

In order to write our $q$-deformed equations in compact form it is
necessary to introduce some additional operators. We first define
the operators: \eqnn{opes} \hm^\pm_{\k}\ \vf \ =& \
\sum_{j,n,\ell,m\in\bbz_+}\ \mu_{j n\ell m}\ \hm^\pm_{\k}\ \vf_{j
n\ell m} \ , \qquad \k = \pm, v,\bv \ ,
\\ [2mm]
T^\pm_{\k}\ \vf \ =& \ \sum_{j,n,\ell,m\in\bbz_+}\ \mu_{j n\ell m}\
T^\pm_{\k}\ \vf_{j n\ell m} \ , \qquad \k = \pm, v,\bv \ ,\nn \eea
and \ $\hm^\pm_{+}\,$, \ $\hm^\pm_{-}\,$, \ $\hm^\pm_{v}\,$,
 \ $\hm^\pm_{\bv}\,$, resp., acts on \ $\vf_{j n\ell m}$ \ by changing
by \ $\pm1$ \ the value of \ $j,n,\ell,m$, resp., while
 \ $T^\pm_{+}\,$, \ $T^\pm_{-}\,$, \ $T^\pm_{v}\,$,
 \ $T^\pm_{\bv}\,$, resp., acts on \ $\vf_{j n\ell m}$ \
by multiplication by \ $q^{\pm j},q^{\pm n},q^{\pm \ell},q^{\pm m}$,
resp. We shall use also the 'logs' $N_\k$ such that ~$T_\k ~=~
q^{N_\k}$. Now we can define the $q$-difference operators:
\eqn{qdif} \hd_{\k} \ \vf \ = \ {1\over \l}\ \hm^{-1}_{\k}\ \left(
T_{\k} - T^{-1}_{\k} \right) \ \vf \  = \ {1\over \l}\
\hm^{-1}_{\k}\ \left( q^{N_{\k}} - q^{-N_{\k}} \right) \ \vf \ . \ee
Note that when $q\ra 1$ then $\hd_{\k}\ra \pd_k\,$. Using
(\ref{opes}) and (\ref{qdif}) the $q$-d'Alembert equation may be
written as \cite{VKD3}, \cite{DGPZ}, respectively, \eqn{qdal} \left(
q\ \hd_-\ \hd_+\ T_v\ T_\bv \ -\ \hd_v\ \hd_\bv \right)\ T_v\ T_-\
T_+\ T_\bv \ \hvf \ \ = \ \ 0 \,, \ee \eqn{qqdal} \left( \hd_-\
\hd_+\ -\ q\ \hd_v\ \hd_\bv \ T_v\ T_\bv \right)\ T_-\ T_+\ \ \tvf \
\ = \ \ 0 \,. \ee Note that when $q\ra 1$ both equations
(\ref{qdal}), (\ref{qqdal}) go to (\ref{dali}). Note that the
operators in (\ref{opes}), (\ref{qdif}), (\ref{qdal}), (\ref{qqdal})
for different variables commute, i.e., we have passed to commuting
variables. However, keeping the normal ordering it is
straightforward to pass back to noncommuting variables.

 Using results from  \cite{VKD4} we have for the $q$-analogue of (\ref{opsb}):
\eqnn{opsc} _qI^+(n) ~&=&~ \ha \Bigl( [n]_q\,[n-1]_q\ _qI^2_1
\,_qI^2_2 - [2]_q\,[n-1]_q\,[n+1]_q \,_qI_1 \,_qI^2_2 \,_qI_1 +\cr
&&+ [n]_q\,[n+1]_q \,_qI^2_2 \,_qI^2_1 \Bigr)
~, \\
_qI^-(n) ~&=&~ \ha \Bigl( [n]_q\,[n-1]_q \,_qI^2_3 \,_qI^2_2 -
[2]_q\,[n-1]_q\,[n+1]_q \,_qI_3 \,_qI^2_2 \,_qI_3  + \cr&&+
[n]_q\,[n+1]_q \,_qI^2_2 \,_qI^2_3 \Bigr) ~,\nn \eea where the
$q$-deformed versions $\,_qI_a$ of (\ref{opsbb}) in the basis
(\ref{spca})  are: \eqnn{rek}
_qI_1\ &=&\ \hd_zT_zT_vT_+(T_-T_{\bv})^{-1} \\
_qI_2\ &=&\ (q\hm_z\hd_vT_-^2+\hm_z\hm_{\bz}\hd_+T_-T_{\bv}T_v^{-1}+
\hd_-T_-\ +\cr &&+\
q^{-1}\hm_{\bz}\hd_{\bv}-\l\hm_v\hm_{\bz}\hd_-\hd_+T_{\bv})\ T_{\bv}
T_{\bz}^{-1} \nn\\
_qI_3\ &=&\ \hd_{\bz}T_{\bz}\ .\nn \eea (For comparison, note that
in  the $q$-Maxwell operators  are used the following expressions:
$\,_qI^+_n = \ha ( [n+2]_q\, _qI_1\, _qI_2  - [n+3]_q\, _qI_2\,
_qI_1 )$, $\,_qI^-_n ~=~ \ha ( [n+2]_q\, _qI_3\, _qI_2  - [n+3]_q\,
_qI_2\, _qI_3 )$.)

Then the $q$-Weyl equations are (cf. (\ref{wgee})): \eqn{wgqq}
_qI^+(4) ~ C^+(z) ~=~ T(z,\bz) ~, \qquad _qI^-(4) ~ C^-(\bz) ~=~
T(z,\bz) ~, \ee  while $q$-analogues of (\ref{wme}) are: \eqn{wmeq}
_qI^+(2) ~ h(z,\bz) ~=~ C^-(\bz) ~, \qquad _qI^-(2) ~ h(z,\bz) ~=~
C^+(z) ~. \ee

We shall look for solutions of the $q$-Weyl gravity equations in
terms of a deformation of the plane wave given in
 \cite{DPZa}. This deformation is given in the basis (\ref{spca}):
\bea\label{dploska} \widehat{\exp}_q (k, x) \ &=&\
\sum_{s=0}^\infty \, {1 \over {[s]_q!}}\ \hh_s \ ,\\
&& [s]_q! \equiv [s]_q [s-1]_q \cdots [1]_q \ , \quad [0]_q! \equiv
1 \ , \quad [n]_q \equiv { q^n-q^{-n}\over q-q^{-1} } \ , \nn\eea
\bea\label{hfs} \hh_s\ &=&\ {\footnotesize
 \beta^s \sum_{a,b,n \in \bbz_+}
{ (-1)^{s-a-b}\ q^{n(s-2a-2b+2n)\ +\ a(s-a-1)\ +\
b(-s+a+b+1)\ +\ P_s(a,b)}\
\over {\G_q(a-n+1)\G_q(b-n+1)\G_q(s-a-b+n+1)[n]_q! }} \ \times }\nn\\
&&\times \ k_v^{s-a-b+n} k_-^{b-n} k_+^{a-n} k_\bv^n v^n x_-
^{a-n} x_+^{b-n} \bv^{s-a-b+n} \ , \\[2mm]
&&\left(\beta^s \right)^{-1} \ = \ \sum_{p=0}^s \ {q^{(s-p)(p-1)+p}
\over [p]_q!\ [s-p]_q!} \ ,\nn \eea where the momentum components \
$(k_v,k_-,k_+,k_\bv)$\ are supposed to be non-commu\-ta\-tive
between themselves (obeying the same rules (\ref{coop}) as the
$q$-Min\-kow\-ski coordinates), and commutative with the
coordinates. Further, $\G_q$ is a $q$-deformation of the
$\G$-function, of which here we use only the properties: $\G_q(p) =
[p-1]_q!$ for $p\in\bbn$, $1/\G_q(p) = 0$ for $p\in\bbz_-\,$;\
$P_s(a,b)$ is a polynomial in $a,b$. Note that $\
(\hh_s)\vert_{q=1}\ =\ (k\cdot x)^s\ $ and thus $\ (\widehat{\exp}_q
(k, x))\vert_{q=1}\ =\ \exp (k\cdot x)\ $. This $q$-plane wave has
some properties analogous to the classical one but is not an
exponent or $q$-exponent, cf. \cite{GR}.
This is enabled also by the fact (true also for $q=1$)
that solving the equations may be done in terms of the components
$\ \hh_s\ $. This deformation of the plane wave generalizes
the original one from \cite{DK} which is obtained by setting $\ P_s(a,b) =0\,$.
Each $\hh_s$ satisfies the $q$-d'Alembert equation (\ref{qdal}) on the
momentum $q$-cone: \be\label{zap}\cl^k_q \ \ \equiv \ \ k_- k_+ \ -
\ q^{-1}\, k_v k_{\bv} \ = \ k_+ k_- \ - \ q\, k_v k_{\bv} \ =\ 0 \
. \ee

\section{Solutions of $q$-Weyl gravity}

We shall use the basis (\ref{spca}). The solutions of the first equation in
(\ref{wgqq}) in the homogeneous case ($T=0$) are:
\be\label{rsc} _qC_0^{+}
=\ \sum_{s=0}^\infty\ {1\over [s]_q!}\ \hat{C}^{+}_{s}  \ ,\ee
\bea\label{rscc} \hat{C}^{+}_{s} \ &=&\
\sum^4_{m=0}\hat{\gamma}^{s+}_m \left(\ \prod^{-m+3}_{i=0}(k_+ -
q^{i+ B_s +s+4}k_{\overline{v}}z)\right)\left(\ \prod^3_{j=-m+4}(k_v
-q^{j+ B_s +s+4}k_-z)\right) \hat{h}^{+}_s \ ,
\eea
where $\ \hh^+_s\ $ is $\ \hh_s\ $ with:
\eqn{ppp}
P_s(a,b) ~=~ P^+_s(a,b) ~\equiv~ R_s(a) + B_s b , \ee
$\hgamma^{s+}_{m}\,, B_s$ are arbitrary constants,
$R_s(a)$ is an arbitrary polynomial in $a$.
Note that the factors preceding $\hh_s^+$ depend
on $B_s$ but not on  $R_s(a)$.
The check that (\ref{rsc}) is a solution is
done for commutative Minkowski coordinates and noncommutative
momenta on the q-cone.
In order to be able to write the above solution
in terms of the deformed plane wave we have to suppose that
the $\hgamma^{s+}_{m}\,,B_s+s$ for different $s$ coincide:
\ $\hgamma^{s+}_{m}=\tgamma^+_{m}\,$, e.g., we can make the
choice $B_s=B'-s-4$. Then we have:
\be\label{teshp} _q C^+_0\ =\ \sum_{m=0}^4\ \tgamma^+_{m}\
\left(\ \prod_{i=0}^{-m+3}(k_+-q^{i+B'} k_\bv z) \right)\,
\left(\ \prod_{j=-m+4}^3 (k_v - q^{j+B'} k_- z) \right)\
\widehat{\exp}_q^+ (k, x)\ ,\ee
where $\ \widehat{\exp}_q^+ (k, x)\ $ is $\ \widehat{\exp}_q (k,
x)\ $ with the choice (\ref{ppp}).

The solutions of the second equation in (\ref{wgqq})  are:
\eqn{sashp} _q C^-_0\ =\ \sum_{s=0}^\infty\ {1\over [s]_q!}\ \hat C^-_s
\ee
\eqnn{rashp} &&\hat C^-_s\ =\ \sum_{m=0}^4\ \hgamma^{s-}_{m}\
\left(\ \prod_{i=-1}^{-m+2}(k_+-q^{i-D_s} k_v \bz)\right)\,
\left(\ \prod_{j=-m+3}^2 (k_\bv - q^{j-D_s} k_-\, \bz) \right)\
\hh_s^- \eea
where $\ \hh_s^-\ $ is $\ \hh_s\ $ with:
\be P_s(a,b) ~=~ P^-_s(a,b) ~\equiv~
D_s a + Q_s(b), \label{cp}\ee
$\hgamma^{s-}_{m}\,,D_s$ are arbitrary constants, and
$Q_s(b)$ is an arbitrary polynomial.
In order to be able to write this solution
in terms of the deformed plane wave we have to suppose that
the $\hgamma^{s-}_{m}\,,D_s$ for different $s$ coincide:
\ $\hgamma^{s-}_{m}=\hgamma^-_{m}\,$, $D_s=D$. Then we have:
\be\label{tashp} _qC^-_0 \ =\ \sum_{m=0}^4\ \hgamma^-_{m}\
\left( \prod_{i=-1}^{-m+2}(k_+-q^{i-D} k_v \bz)\right)\,
\left( \prod_{j=-m+3}^2 (k_\bv - q^{j-D} k_-\, \bz) \right)\
\widehat{\exp}_q^- (k, x) \ , \ee
where $\ \widehat{\exp}_q^- (k, x)\ $ is $\widehat{\exp}_q (k,
x)$ with the choice (\ref{cp}).

\section{Summary and Outlook}

In the present paper we have constructed  $q$-plane wave solutions
of the earlier proposed quantum conformal deformation (\ref{wgqq})
of the equations of linear conformal  Weyl gravity.  We have
restricted ourselves to the vacuum case with energy-momentum tensor
$T=0$. In our further research we plan to consider solutions of
(\ref{wgqq}) in the case $T\neq 0$. Later on we shall look for
solutions of equations (\ref{wmeq}) and explore their consistency
with the solutions of (\ref{wgqq}).

\section*{Acknowledgements}

The authors were supported in part by the Bulgarian National Council
for Scientific Research, grant F-1205/02, and by the European RTN
'Forces-Universe', contract MRTN-CT-2004-005104. The first author
was supported in part also by the Alexander von Humboldt Foundation
in the framework of the Clausthal-Leipzig-Sofia Cooperation.

\end{document}